\documentclass[12pt]{article}
\usepackage[utf8]{inputenc}
\usepackage[T2A]{fontenc}

\usepackage[english]{babel}
\usepackage{fullpage}

\usepackage{cmap}					
\usepackage{mathtext}

\usepackage{amsmath,amsfonts,amssymb,amsthm,mathtools} 
\usepackage{icomma}

\usepackage{mathrsfs}

\usepackage{bbm}
\usepackage{xcolor}
\usepackage{hyperref}
\usepackage{bm}
\usepackage{cite}

\usepackage{algorithm}
\usepackage[noend]{algorithmic}

\usepackage{kamzolov}

\theoremstyle{definition}

\theoremstyle{plain}
\newtheorem{theorem}{Theorem}
\newtheorem{lemma}{Lemma}

\theoremstyle{remark}

\newcommand{\N}{\mathcal{N}}

\title{Near-Optimal Hyperfast Second-Order Method for Convex Optimization and its Sliding}

\newcommand{\beq}{\begin{equation}}
\newcommand{\eeq}{\end{equation}}

\author{
Dmitry Kamzolov \textsuperscript{a} and
Alexander Gasnikov\textsuperscript{a,b,c}.\\
\textsuperscript{a}Moscow Institute of Physics and Technology, Moscow, Russia;\\
\textsuperscript{b} National Research University Higher School of Economics, Moscow, Russia;\\ 
\textsuperscript{c}Institute for Information Transmission Problems RAS, Moscow, Russia. 
}

\begin{document}

\maketitle

\begin{abstract}
  In this paper, we present a new Hyperfast Second-Order Method with convergence rate $O(N^{-5})$ up to a logarithmic factor for the convex function with Lipshitz the third derivative. This method based on two ideas. The first comes from the superfast second-order scheme of Yu. Nesterov (CORE Discussion Paper 2020/07, 2020). It allows implementing the third-order scheme by solving subproblem using only the second-order oracle. This method converges with rate $O(N^{-4})$. The second idea comes from the work of Kamzolov et al. (arXiv:2002.01004). It is the inexact near-optimal third-order method. In this work, we improve its convergence and merge it with the scheme of solving subproblem using only the second-order oracle. As a result, we get convergence rate $O(N^{-5})$ up to a logarithmic factor. This convergence rate is near-optimal and the best known up to this moment.
Further, we investigate the situation when there is a sum of two functions and improve the sliding framework from Kamzolov et al. (arXiv:2002.01004) for the second-order methods.
\end{abstract}

\section{Introduction}
In recent years, it has been actively developing higher-order or tensor methods for convex optimization problems. The primary impulse was the work of Yu. Nesterov \cite{nesterov2019implementable} about the possibility of the implementation tensor method. He proposed a smart regularization of Taylor approximation that makes subproblem convex and hence implementable. Also Yu. Nesterov proposed accelerated tensor methods  \cite{nesterov2018lectures, nesterov2019implementable}, later A. Gasnikov et al.
\cite{bubeck2018near, gasnikov2019optimal, gasnikov2019near, jiang2018optimal}  proposed the near-optimal tensor method via the Monteiro--Svaiter envelope \cite{monteiro2013accelerated} with line-search and got a near-optimal convergence rate up to a logarithmic factor. Starting from 2018-2019 the interest in this topic rises. There are a lot of developments in tensor methods, like tensor methods for Holder continuous higher-order derivatives \cite{grapiglia2019tensor, song2019towards}, proximal methods\cite{doikov2019contracting}, tensor methods for minimizing the gradient norm of convex function \cite{dvurechensky2019near, grapiglia2019tensor}, inexact tensor methods \cite{grapiglia2019inexact, nesterov2020inexact, kamzolov2020optimal}, and near-optimal composition of tensor methods for sum of two functions \cite{kamzolov2020optimal}. There are some results about local convergence and convergence for strongly convex functions \cite{doikov2019local, gasnikov2019optimal, gasnikov2017universal}. See \cite{gasnikov2017universal} for more references on applications of tensor method. 

At the very beginning of 2020, Yurii Nesterov proposed a Superfast Second-Order Method \cite{nesterov2020superfast} that converges with the rate $O(N^{-4})$ for a convex function with Lipshitz third-order derivative. This method uses only second-order information during the iteration, but assume additional smoothness via Lipshitz third-order derivative.\footnote{Note, that for the first-order methods in non-convex case earlier (see, \cite{carmon1711lower} and references therein) it was shown that additional smoothness assumptions lead to an additional acceleration. In convex case, as far as we know these works of Yu. Nesterov \cite{nesterov2020inexact, nesterov2020superfast} are the first ones where such an idea was developed.} Here we should note that for the first-order methods, the worst-case example can't be improved by additional smoothness because it is a specific quadratic function that has all high-order derivatives bounded \cite{nesterov2020inexact}.\footnote{However, there are some results \cite{wilson2019accelerating} that allow to use tensor acceleration for the first-order schemes. This additional acceleration requires additional assumptions on smoothness. More restrictive ones than limitations of high-order derivatives.} But for the second-order methods, one can see that the worst-case example does not have Lipshitz third-order derivative. This means that under the additional assumption, classical lower bound $O(N^{-2/7})$ can be beaten, and Nesterov proposes such a method that converges with $O(N^{-4})$ up to a logarithmic factor. The main idea of this method to run the third-order method with an inexact solution of the Taylor approximation subproblem by method from Nesterov with inexact gradients that converges with the linear speed. By inexact gradients, it becomes possible to replace the direct computation of the third derivative by the inexact model that uses only the first-order information. Note that for non-convex problems previously was proved that the additional smoothness might speed up algorithms \cite{birgin2017worst, bubeck2019complexity, grapiglia2019inexact, nesterov2006cubic, wang2018cubic}. 

In this paper, we propose a Hyperfast Second-Order Method for a convex function with Lipshitz third-order derivative with the convergence rate $O(N^{-5})$ up to a logarithmic factor. For that reason, firstly, we introduce Inexact Near-optimal Accelerated Tensor Method, based on methods from \cite{bubeck2018near, kamzolov2020optimal} and prove its convergence. Next, we apply Bregman-Distance Gradient Method from \cite{grapiglia2019inexact, nesterov2020superfast} to solve Taylor approximation subproblem up to the desired accuracy. This leads us to Hyperfast Second-Order Method and we prove its convergence rate. 
This method have near-optimal convergence rates for a convex function with Lipshitz third-order derivative and the best known up to this moment.

Also we propose a Hyperfast Second-Order Sliding for a sum of two convex functions with Lipshitz third-order derivative. It is based on ideas \cite{lan2016gradient, lan2016accelerated} and \cite{kamzolov2020optimal}. By sliding frameworks we can separate oracle complexities for different functions, it means that we compute only necessary number of derivatives for each functions, like they are separate and not in the sum. So we use sliding for the third-order methods from \cite{kamzolov2020optimal} and solve inner subproblem by Bregman-Distance Gradient Method from \cite{nesterov2020superfast}. As a result, we get method with separate oracle complexity with the convergence rate $O(N^{-5})$ up to a logarithmic factor. This method have near-optimal oracle complexities for a sum of two convex functions with Lipshitz third-order derivative and the best known on that moment.

The paper is organized as follows. In Section 2 we propose Inexact Near-optimal Accelerated Tensor Method and prove its convergence rate. In Section 3 we propose Hyperfast Second-Order Method and get its convergence speed. In Section 4 we assume, that our problem is the sum of two functions. We propose Hyperfast Second-Order Sliding to separate oracle complexities for two functions. In Section 5 we introduce Machine Learning applications. Section 6 is a conclusion.

\section{Inexact Near-optimal Accelerated Tensor Method}

In what follows, we work in a finite-dimensional linear vector space $E=\R^n$, equipped with a Euclidian norm $\|\,\cdot\,\|=\|\,\cdot\,\|_2$.

We consider the following convex optimization problem:
\begin{equation}
\min\limits_{x} f(x),
\label{eq_pr}
\end{equation} 
where $f(x)$ is a convex function with Lipschitz $p$-th derivative, it means that
\begin{equation}
    \|D^p f(x)- D^p f(y)\|\leq L_{p}\|x-y\|.
    \label{def_lipshitz}
\end{equation}
Then Taylor approximation of function $f(x)$ can be written as follows:
\begin{equation*}
    \Omega_{p}(f,x;y)=f(x)+\sum_{k=1}^{p}\frac{1}{k!}D^{k}f(x)\left[ y-x \right]^k, \, y\in \R^n.
\end{equation*}
By \eqref{def_lipshitz} and the standard integration we can get next two inequalities
\begin{equation*}
    |f(y)-\Omega_{p}(f,x;y)|\leq \frac{L_{p}}{(p+1)!}\|y-x\|^{p+1},
\end{equation*}
\begin{equation}
    \|\nabla f(y)- \nabla \Omega_{p}(f,x;y)\|\leq \frac{L_{p}}{p!}\|y-x\|^{p}.
    \label{eq_model_grad}
\end{equation}

Problem \eqref{eq_pr} can be solved by tensor methods \cite{nesterov2019implementable} or its accelerated versions \cite{nesterov2018lectures, gasnikov2019near, bubeck2018near, jiang2018optimal}. 
These  methods have next basic step:
\begin{equation*}
    \label{RCTMforFG}
T_{H_p}(x) = \argmin_{y} \lb \tilde{\Omega}_{p,H_p}(f,x;y) \rb ,
\end{equation*}
where
\begin{equation}
\label{model}
     \tilde{\Omega}_{p,H_p}(f,x;y) = \Omega_{p}(f,x;y) + \frac{H_p}{p!}\|y - x \|^{p+1}.
\end{equation}
For $H_p\geq L_p$ this subproblem is convex and hence implementable. 

But what if we can not solve exactly this subproblem. In paper \cite{nesterov2020superfast} it was introduced Inexact $p$th-Order Basic Tensor Method (BTMI$_p$) and Inexact $p$th-Order Accelerated Tensor Method (ATMI$_p$). They have next convergence rates $O(k^{-p})$ and $O(k^{-(p+1)})$, respectively. In this section, we introduce Inexact $p$th-Order Near-optimal Accelerated Tensor Method (NATMI${_p}$) with improved convergence rate $\tilde{O}(k^{-\frac{3p+1}{2}})$, where $\tilde{O}(\cdot)$ means up to logarithmic factor. It is an improvement of Accelerated Taylor Descent from \cite{bubeck2018near} and generalization of Inexact Accelerated Taylor Descent from \cite{kamzolov2020optimal}. 

Firstly, we introduce the definition of the inexact subproblem solution. 
Any point from the set
\begin{equation}
\label{inexact}
\N_{p,H_p}^{\gamma}(x) = \lb T \in \R^n  \, : \, \|\nabla \tilde{\Omega}_{p,H_p}(f,x;T) \|\leq \gamma \|\nabla f(T)\| \rb
\end{equation}
is the inexact subproblem solution, where $\gamma\in [0; 1]$ is an accuracy parameter. $N_{p,H_p}^{0}$ is the exact solution of the subproblem. 

Next we propose Algorithm \ref{alg_NATMI}.

\begin{algorithm} [H]
\caption{Inexact $p$th-Order Near-optimal Accelerated Tensor Method (NATMI)}\label{alg_NATMI}
	\begin{algorithmic}[1]
		\STATE \textbf{Input:} convex function $f : \R^n \rightarrow \R$ such that $\nabla^p f$ is $L_p$-Lipschitz, $H_p=\xi L_p$ where $\xi$ is a scaling parameter, $\gamma$ is a desired accuracy of the subproblem solution.
		\STATE Set $A_0 = 0, x_0 = y_0$
		\FOR{ $k = 0$ \TO $k = K- 1$ }
		\STATE Compute a pair $\lambda_{k+1} > 0$ and $y_{k+1}\in \R^n$ such that
		\[
\frac{1}{2} \leq \lambda_{k+1} \frac{H_{p} \cdot \|y_{k+1} - \tilde{x}_k\|^{p-1}}{(p-1)!}  \leq \frac{p}{p+1} \,,
\]
where
\begin{equation}
\label{prox_step}
y_{k+1} \in \N_{p,H_p}^{\gamma}(\tilde{x}_k)
\end{equation}
and
		\[
a_{k+1} = \frac{\lambda_{k+1}+\sqrt{\lambda_{k+1}^2+4\lambda_{k+1}A_k}}{2} 
\text{ , } 
A_{k+1} = A_k+a_{k+1}
\text{ , and } 
\tilde{x}_k = \frac{A_k}{A_{k + 1}}y_k + \frac{a_{k+1}}{A_{k+1}} x_k \,.
		\]
		\STATE Update $x_{k+1} := x_k-a_{k+1} \nabla f(y_{k+1})$
		\ENDFOR
		\RETURN $y_{K}$ 
	\end{algorithmic}	
\end{algorithm}

To get the convergence rate of Algorithm \ref{alg_NATMI} we prove additional lemmas. The first lemma gets intermediate inequality to connect theory about inexactness and method's theory.

\begin{lemma}
If $y_{k+1} \in \N_{p,H_p}^{\gamma}(\tilde{x}_k) $, then
\begin{equation}
\label{inexact_norms}
    \|\nabla  \tilde{\Omega}_{p,H_p}(f,\tilde{x}_k;y_{k+1}) \|  \leq \frac{\gamma}{1-\gamma} \cdot \frac{(p+1)H_p+L_p}{p!}\|y_{k+1}-\tilde{x}_k\|^p. 
\end{equation}
\end{lemma}
\begin{proof}
From triangle inequality we get
\begin{align*}
    \| \nabla f(y_{k+1}) \| &\leq \| \nabla f(y_{k+1}) - \nabla \Omega_{p}(f,\tilde{x}_k;y_{k+1}) \|\\
    &+ \|\nabla \Omega_{p}(f,\tilde{x}_k;y_{k+1})-\nabla  \tilde{\Omega}_{p,H_p}(f,\tilde{x}_k;y_{k+1}) \| + \|\nabla  \tilde{\Omega}_{p,H_p}(f,\tilde{x}_k;y_{k+1}) \|\\
    &\overset{\eqref{eq_model_grad},\eqref{model},\eqref{inexact}}{\leq} \frac{L_p}{p!}\| y_{k+1}-\tilde{x}_k \|^{p}+ \frac{(p+1)H_p}{p!}\| y_{k+1} - \tilde{x}_k \|^{p} + \gamma \|\nabla  f(y_{k+1}) \|.
\end{align*}
Hence,
\begin{equation*}
    (1-\gamma)\| \nabla f(y_{k+1}) \| 
    \leq \frac{(p+1)H_p+L_p}{p!}\| y_{k+1} - \tilde{x}_k \|^{p}.
\end{equation*}
And finally from \eqref{inexact} we get
\begin{equation*}
    \|\nabla  \tilde{\Omega}_{p,H_p}(f,\tilde{x}_k;y_{k+1}) \| \leq  \frac{\gamma}{1-\gamma} \cdot \frac{(p+1)H_p+L_p}{p!}\| y_{k+1} - \tilde{x}_k \|^{p}.
\end{equation*}
\end{proof}

Next lemma plays the crucial role in the prove of the Algorithm \ref{alg_NATMI} convergence.
It is the generalization for inexact subpropblem of Lemma 3.1 from \cite{bubeck2018near}.
\begin{lemma} \label{lem:controlstepsize1}
If $y_{k+1} \in \N_{p,H_p}^{\gamma}(\tilde{x}_k) $, $H_p=\xi L_p$ such that $1 \geq  2\gamma+\frac{1}{\xi(p+1)}$ and
\begin{equation} \label{eq_lemma3_HP}
\frac{1}{2} \leq \lambda_{k+1} \frac{H_p \cdot \|y_{k+1} - \tilde{x}_k\|^{p-1}}{(p-1)!}  \leq \frac{p}{p+1} \, ,
\end{equation}
then 
\begin{equation*}
    \|y_{k+1} - (\tilde{x}_k - \lambda_{k+1} \nabla f(y_{k+1})) \| \leq \sigma \cdot \|y_{k+1}-\tilde{x}_k\|
\end{equation*}
and
\begin{equation}
\label{eq_inexact_sigma}
    \sigma \geq \frac{p \xi + 1 -\xi +2\gamma\xi }{(1-\gamma)2p \xi},
\end{equation}
where $\sigma \leq 1$.
\end{lemma}
\begin{proof}
Note, that by definition
\begin{equation}
    \label{eq_lemma3_grad}
\nabla \tilde{\Omega}_{p,H_p}(f,\tilde{x}_k;y_{k+1}) = \nabla \Omega_{p}(f,\tilde{x}_k;y_{k+1}) + \frac{H_p(p+1)}{p!}\|y_{k+1} - \tilde{x}_k \|^{p-1} (y_{k+1}-\tilde{x}_k).
\end{equation}
Hence,
\begin{equation}
    \label{eq_lemma3_y-x}
    y_{k+1}-\tilde{x}_k= \frac{p!}{H_p(p+1)\|y_{k+1} - \tilde{x}_k \|^{p-1}}\ls \nabla \tilde{\Omega}_{p,H_p}(f,\tilde{x}_k;y_{k+1}) - \nabla \Omega_{p}(f,\tilde{x}_k;y_{k+1})\rs.
\end{equation}
Then, by triangle inequality we get
\begin{align*}
& \|y_{k+1} - (\tilde{x}_k - \lambda_{k+1} \nabla f(y_{k+1})) \|  = \| \lambda_{k+1} (\nabla f(y_{k+1})- \nabla \Omega_{p}(f,\tilde{x}_k;y_{k+1}))\\
&+\lambda_{k+1}\nabla \tilde{\Omega}_{p,H_p}(f,\tilde{x}_k;y_{k+1})+ \ls y_{k+1} - \tilde{x}_k + \lambda_{k+1}(\nabla \Omega_{p}(f,\tilde{x}_k;y_{k+1})-\nabla \tilde{\Omega}_{p,H_p}(f,\tilde{x}_k;y_{k+1}))\rs \| \\
& \overset{\eqref{eq_model_grad},\eqref{eq_lemma3_y-x}}{\leq} \lambda_{k+1} \frac{L_p}{p!} \|y_{k+1} - \tilde{x}_k\|^p +  \lambda_{k+1}\|\nabla \tilde{\Omega}_{p,H_p}(f,\tilde{x}_k;y_{k+1})\|\\ 
&+ \left|\lambda_{k+1} - \frac{p!}{H_p \cdot (p+1) \cdot \|y_{k+1} - \tilde{x}_k\|^{p-1}} \right| \cdot \|\nabla \tilde{\Omega}_{p,H_p}(f,\tilde{x}_k;y_{k+1})-\nabla \Omega_{p}(f,\tilde{x}_k;y_{k+1})\|  \\
& \overset{\eqref{inexact_norms},\eqref{eq_lemma3_grad}}{\leq} 
\|y_{k+1} - \tilde{x}_k\|\ls \lambda_{k+1} \frac{L_p}{p!} \|y_{k+1} - \tilde{x}_k\|^{p-1} +  \lambda_{k+1}\frac{\gamma}{1-\gamma} \cdot \frac{(p+1)H_p+L_p}{p!}\|y_{k+1}-\tilde{x}_k\|^{p-1}\rs\\ 
&+ \left|\lambda_{k+1} - \frac{p!}{H_p \cdot (p+1) \cdot \|y_{k+1} - \tilde{x}_k\|^{p-1}} \right| \cdot \frac{(p+1)H_p}{p!} \|y_{k+1} - \tilde{x}_k\|^{p}  \\
&=\|y_{k+1} - \tilde{x}_k\|\ls \frac{\lambda_{k+1}}{p!} \ls L_p +  \frac{\gamma}{1-\gamma} ((p+1)H_p+L_p) \rs  \|y_{k+1} - \tilde{x}_k\|^{p-1} \rs\\ 
&+ \|y_{k+1} - \tilde{x}_k\|\left|\frac{\lambda_{k+1}(p+1)H_p}{p!} \|y_{k+1} - \tilde{x}_k\|^{p-1} - 1\right| \\
&\overset{\eqref{eq_lemma3_HP}}{\leq}\|y_{k+1} - \tilde{x}_k\|\ls \frac{\lambda_{k+1}}{p!} \ls L_p +  \frac{\gamma}{1-\gamma} ((p+1)H_p+L_p) \rs  \|y_{k+1} - \tilde{x}_k\|^{p-1} \rs\\ 
&+ \|y_{k+1} - \tilde{x}_k\|\ls 1-\frac{\lambda_{k+1}(p+1)H_p}{p!} \|y_{k+1} - \tilde{x}_k\|^{p-1} \rs \\
&=\|y_{k+1} - \tilde{x}_k\|\ls 1 + \frac{\lambda_{k+1}}{p!} \ls L_p - (p+1)H_p +  \frac{\gamma}{1-\gamma} ((p+1)H_p+L_p) \rs  \|y_{k+1} - \tilde{x}_k\|^{p-1} \rs.
\end{align*}
Hence, by \eqref{eq_lemma3_HP} and simple calculations we get
\begin{align*}
    \sigma &\geq  1 + \frac{1}{2p H_p} \ls L_p - (p+1)H_p +  \frac{\gamma}{1-\gamma} ((p+1)H_p+L_p) \rs \\
    &=  1 + \frac{1}{2p \xi} \ls 1 - (p+1)\xi +  \frac{\gamma}{1-\gamma} ((p+1)\xi +1) \rs \\
    &=  1 + \frac{1}{2p \xi} \ls 1 - p\xi-\xi +  \frac{\gamma p\xi+\gamma \xi +\gamma}{1-\gamma}  \rs \\
    &=  1 + \frac{1}{2p \xi} \ls \frac{1 - p\xi-\xi - \gamma + \gamma p\xi+\gamma\xi+\gamma p\xi+\gamma \xi +\gamma}{1-\gamma}  \rs \\
    &=  1 +  \ls \frac{1 - p\xi-\xi + 2\gamma p\xi+2\gamma\xi }{(1-\gamma)2p \xi}  \rs \\
    &=  \frac{p \xi + 1 -\xi +2\gamma\xi }{(1-\gamma)2p \xi}.  
\end{align*}
Lastly, we prove that $\sigma\leq 1$. For that we need
\begin{align*}
    (1-\gamma)2p \xi &\geq  p \xi + 1 -\xi +2\gamma\xi\\ 
    (p+1) \xi &\geq   1  +2\gamma\xi(1+p)\\
    \frac{1}{2}-\frac{1}{2\xi(p+1)} &\geq  \gamma.
\end{align*}
\end{proof}
We have proved the main lemma for the convergence rate theorem, other parts of the proof are the same as \cite{bubeck2018near}. As a result, we get the next theorem.

\begin{theorem} \label{theoremNATMI}
 Let $f$ be a convex function whose $p^{th}$ derivative is $L_p$-Lipschitz and $x_{\ast}$ denote a minimizer of $f$. Then Algorithm \ref{alg_NATMI} converges with rate
 \begin{equation*}
 f(y_k) - f(x_{\ast}) \leq \tilde{O}\ls\frac{H_p R^{p+1}}{k^{\frac{ 3p +1}{2}}}\rs \,,
 \end{equation*}
 where 
\begin{equation*}
R=\|x_0 - x^{\ast}\|
\end{equation*} 
is the maximal radius of the initial set.
\end{theorem}

\section{Hyperfast Second-Order Method}
In recent work \cite{nesterov2020superfast} it was mentioned that for convex optimization problem \eqref{eq_pr} with first order oracle (returns gradient) the well-known complexity bound $\left(L_{1}R^2/\e\right)^{1/2}$ can not be beaten even if we assume that all $L_{p} < \infty$. This is because of the structure of the worth case function
$$f_p(x) = |x_1|^{p+1} + |x_2 - x_1|^{p+1} + ... + |x_n - x_{n-1}|^{p+1},$$
where $p = 1$ for first order method. It's obvious that $f_p(x)$ satisfy the condition $L_{p} < \infty$ for all natural $p$. So additional smoothness assumptions don't allow to accelerate additionally.  The same thing takes place, for example, for $p=3$. In this case, we also have $L_{p} < \infty$ for all natural $p$. But what is about $p=2$? In this case $L_3 = \infty$. It means that $f_2(x)$ couldn't be the proper worth case function for the second order method with additional smoothness assumptions. So there appears the following question: Is it possible to improve the bound $\left(L_{2}R^3/\e\right)^{2/7}$? At the very beginning of 2020 Yu. Nesterov gave a positive answer. For this purpose, he proposed to use an accelerated third-order method that requires  $\tilde{O}\left((L_{3}R^4/\e)^{1/4}\right)$ iterations by using second-order oracle \cite{nesterov2019implementable}. So all this means that if $L_3 < \infty$, then there are methods that can be much faster than $\tilde{O}\left(\left(L_{2}R^3/\e\right)^{2/7}\right)$.

In this section, we improve convergence speed and reach near-optimal speed up to logarithmic factor. We consider problem \eqref{eq_pr} with $p=3$, hence $L_3<\infty$.
In previous section, we have proved that Algorithm \ref{alg_NATMI} converges. Now we fix the parameters for this method
\begin{equation*}
    p=3,\quad \gamma=\frac{1}{2p}=\frac{1}{6}, \quad \xi = \frac{2p}{p+1}=\frac{3}{2}. 
\end{equation*}
By \eqref{eq_inexact_sigma} we get $\sigma = 0.6$ that is rather close to initial exact $\sigma_{0}=0.5$. For such parameters we get next convergence speed of Algorithm \ref{alg_NATMI} to reach accuracy $\e$:
\begin{equation*}
    N_{out}= \tilde{O}\ls\ls \frac{L_3 R^{4}}{\e}\rs^{\frac{1}{5}}\rs.
\end{equation*}

Note, that at every step of Algorithm \ref{alg_NATMI} we need to solve next subproblem with accuracy $\gamma=1/6$
\begin{equation}
\label{hyper_sub}
  \argmin_{y} \lb \la\nabla f(x_i),y-x_i\ra +\frac{1}{2}\nabla^2 f(x_i)[y-x_i]^2+\frac{1}{6}D^3f(x_i)[y-x_i]^3 + \frac{L_3}{4}\|y - x_i \|^{4} \rb.
\end{equation}

In \cite{grapiglia2019inexact} it was proved, that problem \eqref{hyper_sub} can be solved by Bregman-Distance Gradient Method (BDGM) with linear convergence speed. According to \cite{nesterov2020superfast} BDGM can be improved to work with inexact gradients of the functions. This made possible to approximate $D^3 f(x)$ by gradients and escape calculations of $D^3 f(x)$ at each step. As a result, in \cite{nesterov2020superfast} it was proved, that subproblem \eqref{hyper_sub} can be solved up to accuracy $\gamma = 1/6$ with one calculation of Hessian and $O\ls\log \ls\frac{\|\nabla f(x_i)\|+\|\nabla^2 f(x_i)\|}{\e}\rs\rs$ calculation of gradient.

We use BDGM to solve subproblem from Algorithm \ref{alg_NATMI} and, as a result, we get next Hyperfast Second-Order method as merging NATMI and BDGM.
\begin{algorithm} [H]
\caption{Hyperfast Second-Order Method}\label{alg_hyper}
	\begin{algorithmic}[1]
	    \STATE \textbf{Input:} convex function $f : \R^n \rightarrow \R$ with $L_3$-Lipschitz $3$rd-order derivative.
		\STATE Set $A_0 = 0, x_0 = y_0$
		\FOR{ $k = 0$ \TO $k = K- 1$ }
		\STATE Compute a pair $\lambda_{k+1} > 0$ and $y_{k+1}\in \R^n$ such that
		\[
\frac{1}{2} \leq \lambda_{k+1} \frac{3 L_3 \cdot \|y_{k+1} - \tilde{x}_k\|^{2}}{4}  \leq \frac{3}{4} \,,
\]
where $y_{k+1} \in \N_{3,3L_3/2}^{1/6}(\tilde{x}_k)$ solved by Algorithm \ref{alg_BDGM} and
		\[
a_{k+1} = \frac{\lambda_{k+1}+\sqrt{\lambda_{k+1}^2+4\lambda_{k+1}A_k}}{2} 
\text{ , } 
A_{k+1} = A_k+a_{k+1}
\text{ , and } 
\tilde{x}_k = \frac{A_k}{A_{k + 1}}y_k + \frac{a_{k+1}}{A_{k+1}} x_k \,.
		\]
		\STATE Update $x_{k+1} := x_k-a_{k+1} \nabla f(y_{k+1})$
		\ENDFOR
		\RETURN $y_{K}$ 
	\end{algorithmic}	
\end{algorithm}

\begin{algorithm} [H]
\caption{Bregman-Distance Gradient Method}\label{alg_BDGM}
	\begin{algorithmic}[1]		
		\STATE Set $z_0=\tilde{x}_k$ and $\tau=\frac{3\delta}{8(2+\sqrt{2})\|\nabla f(\tilde{x}_k)\|} $
		\STATE Set objective function
		$$
		\varphi_k(z)= \la\nabla f(\tilde{x}_k),z-\tilde{x}_k\ra +\nabla^2 f(\tilde{x}_k)[z-\tilde{x}_k]^2+D^3f(\tilde{x}_k)[z-\tilde{x}_k]^3 + \frac{L_3}{4}\|z - \tilde{x}_k \|^{4} 
		$$
		\STATE Set feasible set
		\begin{equation*}
		S_k = \lb z: \|z- \tilde{x}_k\|\leq 2 \ls \frac{2+\sqrt{2}}{L_3}\|\nabla f(\tilde{x}_k)\|\rs^{\frac{1}{3}} \rb
		\end{equation*}
		\STATE Set scaling function
		\begin{equation*}
		\rho_k(z)= \frac{1}{2}\la\nabla^2 f(\tilde{x}_k)(z-\tilde{x}_k),z-\tilde{x}_k\ra + \frac{L_3}{4} \|z - \tilde{x}_k \|^4 
		\end{equation*}
		\FOR{ $k \geq 0$  }
		\STATE Compute the approximate gradient $g_{\varphi_k,\tau}(z_i)$ by \eqref{eq_g}.
		\STATE \textbf{IF} $\|g_{\varphi_k,\tau}(z_i)\|\leq \frac{1}{6}\|\nabla f(z_i) \| - \delta$, then \textbf{STOP}
		\STATE \textbf{ELSE} $z_{i+1}=\argmin\limits_{z\in S_k} \lb \la g_{\varphi_k,\tau}(z_{i}),z-z_{i}\ra+2\ls 1 + \frac{1}{\sqrt{2}} \rs\beta_{\rho_k}(z_i,z)\rb,$
		\ENDFOR
		\RETURN $z_{i}$ 
	\end{algorithmic}	
\end{algorithm}

Here $\beta_{\rho_k}(z_i,z)$ is a Bregman distance generated by $\rho_k(z)$
\begin{equation*}
    \beta_{\rho_k}(z_i,z)=\rho_k(z) - \rho_k(z_i) -\la \nabla\rho_k(z_i), z-z_i \ra.
\end{equation*}
By $g_{\varphi_k,\tau}(z)$ we take an inexact gradient of the subproblem \eqref{hyper_sub}
\begin{equation}
    \label{eq_g}
    g_{\varphi_k,\tau}(z)= \nabla f(\tilde{x}_k) +\nabla^2 f(\tilde{x}_k)[z-\tilde{x}_k]+g^{\tau}_{\tilde{x}_k}(z) + L_3\|z - \tilde{x}_k \|^{2} (z - \tilde{x}_k)
\end{equation}
and $g^{\tau}_{\tilde{x}_k}(z)$ is a inexact approximation of $D^3f(\tilde{x}_k)[y-\tilde{x}_k]^2$
\begin{equation*}
    g^{\tau}_{\tilde{x}_k}(z)= \frac{1}{\tau^2}\ls \nabla f(\tilde{x}_k+\tau(z-\tilde{x}_k))+ \nabla f(\tilde{x}_k-\tau(z-\tilde{x}_k))-2\nabla f(\tilde{x}_k)\rs.
\end{equation*}
In paper \cite{nesterov2020superfast} it is proved, that we can choose 
$$\delta=O\ls\frac{\e^{\frac{3}{2}}}{\|\nabla f(\tilde{x}_k)\|^{\frac{1}{2}}_{\ast}+\|\nabla^2 f(\tilde{x}_k)\|^{\frac{3}{2}}/L_3^{\frac{1}{2}}} \rs,$$
then total number of inner iterations equal to
\begin{equation*}
    T_k(\delta)=O\ls\ln{\frac{G+H}{\e}}\rs,
\end{equation*}
where $G$ and $H$ are the uniform upper bounds for the norms of the gradients and Hessians computed at the points generated by the main algorithm.
Finally, we get next theorem.
\begin{theorem} 
 Let $f$ be a convex function whose third derivative is $L_3$-Lipschitz and $x_{\ast}$ denote a minimizer of $f$. Then to reach accuracy $\e$ Algorithm \ref{alg_hyper} with Algorithm \ref{alg_BDGM} for solving subproblem computes 
 \begin{equation*} 
 N_{1}=\tilde{O}\ls\ls \frac{L_3 R^{4}}{\e}\rs^{\frac{1}{5}}\rs 
 \end{equation*}
 Hessians and
 \begin{equation*} 
 N_{2}=\tilde{O}\ls\ls \frac{L_3 R^{4}}{\e}\rs^{\frac{1}{5}}\log\ls\frac{G+H}{\e}\rs\rs 
 \end{equation*}
 gradients, where $G$ and $H$ are the uniform upper bounds for the norms of the gradients and Hessians computed at the points generated by the main algorithm.
\end{theorem}

One can generalize this result on uniformly-strongly convex functions by using inverse restart-regularization trick from \cite{gasnikov2018hypothesis}.

So, the main observation of this section is as follows: If $L_3 < \infty$, then we can use this superfast\footnote{Here we use terminology introduced in \cite{nesterov2020superfast}.} second order algorithm instead of considered in the paper optimal one to make our sliding faster (in convex and uniformly convex cases).

\section{Hyperfast Second-order Sliding}
In this section, we consider problem
\begin{equation*}
\min\limits_{x} f(x)= g(x)+h(x),
\end{equation*} 
where $g(x)$ and $h(x)$ are the convex functions such that $\nabla^3 g$ is $L_{3,g}$-Lipschitz and $\nabla^3 h$ is $L_{3,h}$-Lipschitz, also $L_{3,g}\leq L_{3,h}$ .

In \cite{kamzolov2020optimal} it was proposed an algorithmic framework for composition of tensor methods, also called sliding. This framework separates oracle complexity, hence, we get much smaller number calls for function $g(x)$.
In this section we combine sliding technique and Hyperfast Second-Order Method \ref{alg_hyper} to get Hyperfast Second-Order Sliding.

Firstly, we need NATMI version with smooth composite part as an outer basic method.
\begin{algorithm} [H]
\caption{$3$rd-Order NATMI with smooth composite part}\label{alg_slide_out}
	\begin{algorithmic}[1]
	    \STATE \textbf{Input:} convex functions $h(x)$ and $g(x)$ such that $\nabla^3 h$ is $L_{3,h}$-Lipschitz and $\nabla^3 g$ is $L_{3,g}$-Lipschitz.
		\STATE Set $A_0 = 0, x_0 = y_0$
		\FOR{ $k = 0$ \TO $k = K- 1$ }
		\STATE Compute a pair $\lambda_{k+1} > 0$ and $y_{k+1}\in \R^n$ such that
		$$
\frac{1}{2} \leq \lambda_{k+1} \frac{3 L_{3,g} \cdot \|y_{k+1} - \tilde{x}_k\|^{2}}{4}  \leq \frac{3}{4} \,,
$$
where $y_{k+1} \in \lb T \in E  \, : \, \|\nabla \tilde{\Omega}_{3,3L_{3,g}/2}(g,\tilde{x}_k;T)+\nabla h(T) \|\leq \frac{1}{6} \|\nabla f(T)\| \rb$  

solved by Algorithm \ref{alg_BDGM} and
		\[
a_{k+1} = \frac{\lambda_{k+1}+\sqrt{\lambda_{k+1}^2+4\lambda_{k+1}A_k}}{2} 
\text{ , } 
A_{k+1} = A_k+a_{k+1}
\text{ , and } 
\tilde{x}_k = \frac{A_k}{A_{k + 1}}y_k + \frac{a_{k+1}}{A_{k+1}} x_k \,.
		\]
		\STATE Update $x_{k+1} := x_k-a_{k+1} \nabla f(y_{k+1})-a_{k+1} \nabla h(y_{k+1})$
		\ENDFOR
		\RETURN $y_{K}$ 
	\end{algorithmic}	
\end{algorithm}	
The convergence of this method can be easily proved by united proofs of NATMI with CATD from \cite{kamzolov2020optimal}. We get the same the convergence rate. 

Now we propose Hyperfast Second-Order Sliding Method. It contains three levels of methods outer Algorithm \ref{alg_slide_out} for function $g(x)$ and $h(x)$ as a composite part, in the middle Algorithm \ref{alg_slide_out} solves outer's subproblem with model of $g(x)$ and $h(x)$ and the deepest Algorithm \ref{alg_BDGM} solves inner subproblem of the sum of two models of $g(x)$ and $h(x)$.

\begin{algorithm} [H]
\caption{Hyperfast Second-Order Sliding}\label{alg_sliding}
	\begin{algorithmic}[1]
		\STATE \textbf{Input:} convex functions $g(x)$ and $h(x)$ such that $\nabla^3 h$ is $L_{3,h}$-Lipschitz and $\nabla^3 g$ is $L_{3,g}$-Lipschitz.
		\STATE Set $z_0=y_0=x_0$
		\FOR{$k = 0, $ \TO $K-1$}
		\STATE Run Algorithm \ref{alg_slide_out} for problem $g(x)+h(x)$, where $h(x)$ is a composite part.
		        \FOR {$m = 0, $ \TO $M-1$}
		            \STATE Run Algorithm \ref{alg_slide_out} up to desired accuracy for subproblem  \[\min\limits_{y}\ls\tilde{\Omega}_{3,3L_{3,g}/2}(g,\tilde{x}_k;y)+h(y)\rs\]
		        \ENDFOR
		\ENDFOR
		\RETURN $y_{K}$ 
	\end{algorithmic}	
\end{algorithm}

The total complexity is a multiplication of complexities of each submethod in Algorithm \ref{alg_sliding}. Hence, from \cite{kamzolov2020optimal} we get next convergence speed for outer and middle method. To reach $f(x_N) - f(x^{\ast}) \leq \e$, we need $N_g$  computations of derivatives $g(x)$  and $N_h$ computation of derivatives $h(x)$, where
\begin{align}
    \label{convN_g}
    &N_g = \tilde{O}\left[\left( \frac{L_{3,g} R^{4}}{\e} \right)^{\frac{1}{5}}\right],\\
    \label{convN_h}
    &N_h = \tilde{O}\left[\left( \frac{L_{3,h} R^{4}}{\e} \right)^{\frac{1}{5}}\right].
\end{align}
Also from previous section we have known convergence rate of inner Algorithm \eqref{alg_BDGM}, that equal to $O\ls\ln{\frac{G+H}{\e}}\rs$. This lead us to the next theorem.

\begin{theorem}
Assume $f(x)$ and $g(x)$ are convex functions with  Lipshitz third derivative and $L_{3,g}<L_{3,h}$. Then Method \ref{alg_sliding} converges to $f(x_N)-f(x^{\ast})\leq \e$ with $N_g$ as \eqref{convN_g} computations of Hessian $f(x)$,  $N_h$ as \eqref{convN_h} computation of Hessian $h(x)$,  $O\ls N_g \ln{\frac{G+H}{\e}}\rs$ computations of gradients $g(x)$ and $O\ls N_h \ln{\frac{G+H}{\e}}\rs$ computations of gradients $h(x)$, where $G$ and $H$ are the uniform upper bounds for the norms of the gradients and Hessians $f(x)$ computed at the points generated by the main algorithm.
\end{theorem}

One can generalize this result on sum of $n$ functions sorted by $L_{3,f_i}$ and applied consecutively, also it is possible to separate it by some batches.

\section{Machine Learning Applications}
Many modern machine learning applications reduce at the end to the following optimization problem 
\begin{equation*}
    \min_{x\in\R^n} \frac{1}{m}\sum_{k=1}^m f_k(x) + g(x).
\end{equation*}
Sometimes these problems are convex and have uniformly bounded high-order derivatives. For example, in click prediction model\footnote{
\href{http://towardsdatascience.com/mobile-ads-click-through-rate-ctr-prediction-44fdac40c6ff}{http://towardsdatascience.com/mobile-ads-click-through-rate-ctr-prediction}}
we have to solve logistic regression problem $f_k(x) = \log\left(1 + \exp(-y_k\langle a_k, x \rangle)\right)$ with (quadratic) regularizer $g(x)$ \cite{shalev2014understanding}. So the most appropriate conditions to apply the developed second-order scheme are: the complexity of calculating Hessian\footnote{Here $T_{\nabla f_k}$ -- is the complexity of $\nabla f_k(x)$ calculation. But this formula sometimes is just a rough upper bound. For example, if we consider logistic regression and matrix $A = [a_1,...,a_m]$ is sparse we can improve the bound as follows $O\left(ms^2\right)$, where $s$ is upper bound of number of nonzero elements in $\{a_k\}_{k=1}^m$. For click prediction model $s \ll n$.}
$O\left(m\cdot n \cdot T_{\nabla f_k}\right)$ is of the same order as the complexity of inverting Hessian $O(n^{2.5})$ a.o. (arithmetic operations),\footnote{Standard Gaussian method with the complexity $O(n^{3})$ a.o. is not optimal. The best known theoretical method requires $O(n^{2.37})$ a.o. \cite{le2014powers} are not practical one. So we lead the complexity that corresponds to the best theoretical complexity among practical methods. The result $O(n^{2.5})$ a.o. seems to be a folklore result. Some specialists consider that the best practical methods have the complexity close to Gaussian method $O(n^3)$ a.o. For our purposes (see bellow) this is even better if we want to choose $m$ bigger.  Note also that the complexity of auxiliary problem determines not only the complexity of Hessian inversion. There appears also an additional logarithmic factor \cite{nesterov2019implementable}, but we skip it for simplicity.} the number of terms $m$ is not to big,\footnote{Without loss of generality in machine learning applications every time one may consider $m$ to be not large than $\tilde{O}\left(\e^{-2}\right)$ in convex (but non strongly convex) case \cite{feldman2019high, shalev2009stochastic}. Note, that it's possible due to the proper regularization. Note also, that this remark and the sensitivity of tensor methods to gradient estimation \cite{ghadimi2019generalized, wang2018cubic, wang2018note} assume that if $m$ is bigger than $\tilde{O}\left(\e^{-2}\right)$, one can randomly select a subsum with $\tilde{O}\left(\e^{-2}\right)$ terms and solve only this problem instead of using batched gradient, hessian etc.} that is the first-order variance reduction schemes\footnote{As far as we know, for the moment there are no any high-order variance reduction schemes.} \cite{lan2018optimal} $O\left(\left(m + \sqrt{mL_1R^2/\e}\right)T_{\nabla f_k}\right)$  or SGD \cite{duchi2018introductory} $O\left((L_0^2R^2/\e^2) T_{\nabla f_k}\right)$ are don't dominate hyper-fast second order scheme. In particular, for click prediction model if $m = \Omega\left(n^{2.5}/s^2\right)$ the complexity of proposed hyper-fast scheme is\footnote{Quadratic regularizer don't play any role here until it is not strongly convex.} $\tilde{O}\left(\left(ms^2 + n^{2.5}\right)\e^{-1/5}\right) = \tilde{O}\left( ms^2\e^{-1/5}\right)$ that can be better than the complexity of optimal variance reduction scheme $\tilde{O}\left(ms + s\sqrt{m/\e}\right)$
and SGD $O\left(s\e^{-2}\right)$ depends on the relation between $m$, $s$ and $\e$ (here $\e$ is a `relative' accuracy). For example, if $s=O(1)$, $\e = O(m^{-5/3})$ hyper-fast scheme is the best one. Unfortunately, the requirement on accuracy $\e$ is not very practical one.\footnote{However in distributed context by using statistical preconditioned algorithms one may significantly reduce $m$ \cite{hendrikx2020statistically}. So this estimate became practical one.} But, it's important to note, that if $\e = O(m^{-5/9})$ hyper-fast scheme is still better than SGD.\footnote{For SGD we can parallelize calculations on $O\left(\sigma^2R^2/\e^2\large{/}\sqrt{L_1R^2/\e}\right)$ nodes \cite{dvinskikh2019decentralized}, \cite{woodworth2018graph}. But, if we have nonsmooth case ($L_0 < \infty$, but $L_1 = \infty$) or $n$ is small enough, the best known for us strategy is to use ellipsoid method with batched gradient. The number of iterations $\tilde{O}(n^2)$ (it seems that this number can be improved to $\tilde{O}(n)$ for more advanced schemes \cite{lee2019solving}) and at each iteration one should calculate a batch of size $\tilde{O}(\sigma^2R^2/\e^{2})$ .} So if we allow to use parallel calculations variance reduction scheme fails and we should compare our approach only with SGD. In this case we have rather reasonable requirements on accuracy.

Note, that all these formulas can be rewritten in strongly convex case. 

Note also, that in non convex case we hope that the proposed approach sometimes can also works good in practice due to the theoretical results concern auxiliary problem for tensor methods in non convex case  \cite{grapiglia2019inexact} and optimistic practical results concern applications of Monteiro--Svaiter accelerated proximal envelope for non convex problems \cite{ivanova2019adaptive}.

\section{Conclusion}
In this paper, we present Inexact Near-optimal Accelerated Tensor Method and improve its convergence rate. This improvement make it possible to solve the Taylor approximation subproblem by other methods. Next, we propose Hyperfast Second-Order Method and get its convergence speed $O(N^{-5})$ up to logarithmic factor. This method is a combination of Inexact Third-Order Near-Optimal Accelerated Tensor Method with Bregman-Distance Gradient Method for solving inner subproblem.   Finally, we assume, that our problem is the sum of some functions. We propose Hyperfast Second-Order Sliding to separate oracle complexities of two functions. All this methods have near-optimal convergence rates for given problem classes and the best known on that moment.

Though it seems there is no way to obtain super-fast first and zero order schemes we may use the proposed approach to obtain new classes of first and zero order schemes. Namely, for convex smooth enough problems (analogues for strongly convex) it is possible to propose  first-order method that requires $\tilde{O}(n{\e}^{-1/5})$ gradient calculations and that $\tilde{O}({\e}^{-1/5})$ times inverts Hessian (with complexity $O(n^{2.5})$ a.o. of each inversion). This results is not optimal from the modern complexity theory point of view. There exists a theoretical algorithm \cite{lee2019solving} that requires $\tilde{O}(n)$ gradient calculations and $\tilde{O}(n^2)$ a.o. per each iteration. But for the current moment of time this result is far from to be practical one. 

Analogously, one can show that it is possible to propose zero-order method that requires $\tilde{O}(n^2{\e}^{-1/5})$ function value calculations and that $\tilde{O}({\e}^{-1/5})$ times inverts Hessian. The same things about optimality we may say here. But instead of \cite{lee2019solving} one should refers to \cite{lee2019efficient}.

In this paper, we developed near-optimal Hyperfast Second-Order method for sufficiently smooth convex problem in terms of convergence in function. Based on the technique from the work \cite{dvurechensky2019near}, we can also developed  near-optimal Hyperfast Second-Order method for sufficiently smooth convex problem in terms of convergence in the norm of the gradient. In particular, based on the work \cite{guminov2019accelerated} one may show that the complexity of this approach to the dual problem for $1$-entropy regularized optimal transport problem will be $\tilde{O}\left(\left((\sqrt{n})^{4}/\e\right)^{1/5}\right)\cdot O(n^{2.5}) = O(n^{2.9}\e^{-1/5})$ a.o., where $n$ is the linear dimension of the transport plan matrix, that could be better than the complexity of accelerated gradient method and accelerated Sinkhorn algorithm  $O(n^{2.5}\e^{-1/2})$ a.o. \cite{dvurechensky2018computational, guminov2019accelerated}. Note, that the best theoretical bounds for this problem are also far from to be practical ones \cite{blanchet2018towards, jambulapati2019direct, lee2019solving, quanrud2018approximating}. 

Note, that in  March 2020 there appears a new preprint of Yu. Nesterov where the same result $\tilde{O}(\e^{-1/5})$ was obtained based on another proximal accelerated envelop \cite{nesterov2020inexact10}. Note, that both approaches (described in this paper and in \cite{nesterov2020inexact10}) have also the same complexity in terms of $O(\text{ })$ (have factors $\ln^2 \e^{-1}$ due to the line search and requirements to the accuracy we need to solve auxiliary problem).

\subsection{Acknowledgements}
We would like to thank Yurii Nesterov, Pavel Dvurechensky and Cesar Uribe for fruitful discussions. We also would like to thanks Soomin Lee (Yahoo), Erick Ordentlich (Yahoo) and Andrey Vorobyev (Huawei), Evgeny Yanitsky (Huawei), Olga Vasiukova (Huawei) for motivating us by concrete problems formulations.
 
\subsection{Funding}
The work of D. Kamzolov was funded by RFBR, project number 19-31-27001. The work of A.V. Gasnikov in the first part of the paper was supported by RFBR grant 18-29-03071 mk. A.V. Gasnikov was also partially supported by the Yahoo! Research Faculty Engagement Program.

\bibliographystyle{plain} 
\bibliography{biblio}

\end{document}